\documentclass[12pt]{article}
\usepackage{amssymb}
\def \C{{\mathbb C}}
\def \R{{\mathbb R}}
\def \Z{{\mathbb Z}}
\def \p{{\mathbb P}}
\def \Q{{\mathbb Q}}

\begin{document}
\title{Fibrations m\'eromorphes sur certaines vari\'et\'es \`a fibr\'e canonique trivial}
\author{Ekaterina Amerik\thanks{Universit\'e Paris-Sud, Laboratoire des Math\'ematiques,
B\^atiment 425, 91405 Orsay, France. Ekaterina.Amerik@math.u-psud.fr} \ et Fr\'ed\'eric Campana\thanks{Institut Elie Cartan, Universit\'e Nancy I, B.P. 239, 54506
Vandoeuvre-l\`es-Nancy Cedex. campana@iecn.u-nancy.fr}}
\date{}
\maketitle

\ \ \ \ \ \ \ \ \ \ \ \ \ \ \ \ \ \ \ \ \ \ \ \ \ \ \ \ \ \ {\it  A F.A. Bogomolov, en l'honneur de
ses 60 ans}

\section{Introduction}

\par\bigskip
 Cet article pr\'esente quelques observations motiv\'ees
par la 
recherche de vari\'et\'es lisses
projectives complexes $X$ munies d'un endomorphisme $f:X\to X$ de degr\'e $d$ au moins $2$. 

\

Si l'application $f$ est r\'eguli\`ere, il semble difficile de
produire des exemples qui ne soient 
pas d\'eduits des exemples \'evidents (vari\'et\'es toriques ou
ab\'eliennes), 
par produits, passage au quotient par un groupe fini, et cetera. 
Ce qu'indiquent \cite{ARV}, \cite{A}, \cite{F}, \cite{N}.

\

Lorsque $f$ est
rationnelle, il y a de nouveaux exemples, mais certains d'entre eux
sont encore des versions relatives des pr\'ec\'edents. 
Par exemple, si $S$ est une surface projective de dimension canonique, ou de
Kodaira,
\'egale \`a $1$, il est
facile de voir que $S$ admet un endomorphisme rationnel de degr\'e  $d>1$
 (que $S$ ait ou non une section). Un tel endomorphisme $f$ doit
cependant ``respecter'' 
la fibration
elliptique $\phi$ (d'Iitaka-Moishezon), c'est-\`a-dire, envoyer une
fibre de $\phi$
sur une fibre de $\phi$. Si la base est de genre $\geq 2$, une puissance de
$f$
pr\'eserve donc la fibration.

\
 
Si $S$ est une surface
K3 sp\'eciale (par exemple, elliptique ou de
Kummer), elle admet encore de tels endomorphismes rationnels. Si $S$ 
est g\'en\'erique (telle que $Pic(S)=\Z$, disons), on ne conna\^{\i}t,
par contre, aucun exemple (mais il semble difficile de d\'emontrer 
la non-existence).

\

Cependant, il est observ\'e dans \cite{V} qu'il existe
des familles $\{X_t, t\in T\}$ de vari\'et\'es $X_t$
avec $K_{X_t}=0$, et $Pic(X_t)=\Z$ pour $t\in T$ g\'en\'erique,
munies d'un endomorphisme $f_t$ 
de degr\'e $d>1$ pour tout $t\in T$. On peut prendre pour $X_t$ la vari\'et\'e de Fano des $k$-plans
dans une hypersurface cubique lisse $V_t$ de $\p^n$, pour $k$ et $n$ convenables.

\

Par exemple, lorsque $k=1$ et $n=5$, il s'agit de la vari\'et\'e
$X_t={\cal F}(V_t)$ 
(${\cal F}$ pour ``Fano'') des
droites
d'une cubique lisse $V_t$ de dimension 4; $X_t$ est alors lisse de dimension
$4$. De plus, $X_t$ est hyperk\"ahlerienne
irr\'eductible, 
\'equivalente par d\'eformation \`a $Hilb^2(S)$
o\`u
$S$ est une surface K3 de degr\'e 14 (\cite{BD}). Pour $t$
g\'en\'erique, $Pic(X_t)\cong \Z$. L'endomorphisme $f$ associe \`a
$l$, droite g\'en\'erique de 
$V_t\subset \p^5$, la droite $l'$  r\'esiduelle \`a $l$ dans
$P\cap V_t$, o\`u $P$ est l'unique plan dans $\p^5$ tangent \`a $V_t$ le long de
$l$ (on calcule facilement $N_{l,V}\cong {\cal O}_{\p^1}\oplus {\cal
  O}_{\p^1}
\oplus {\cal O}_{\p^1}(1)$, d'o\`u l'existence et l'unicit\'e de $P$) . Dans \cite{V}, Voisin montre que le degr\'e $d$ de $f$ est 16.

\

Une question naturelle est de savoir si ces exemples sont
``nouveaux'' dans le sens de 
la discussion ci-dessus. Et, en effet, nous montrons tout d'abord que, pour $X$ un membre
g\'en\'eral d'une famille comme ci-dessus, un endomorphisme rationnel de $X$
ne pr\'eserve aucune fibration non-triviale.

Plus pr\'ecisement, nous d\'emontrons le r\'esultat suivant:

\medskip

{\bf Th\'eor\`eme 2.1} {\it Soit $X$ une vari\'et\'e lisse projective,
telle que $K_X= 0$ et le groupe de N\'eron-Severi $NS(X)\cong \Z$. Alors toute fibration rationnelle
$g:X \dasharrow B$, $0<dim(B)<dim(X)$, est en vari\'et\'es de type
g\'en\'eral.}

\

Dans la section 4, nous associons une fibration m\'eromorphe 
\`a tout endomorphisme m\'eromophe
dominant:

\medskip

{\bf Th\'eor\`eme 4.1} {\it Soit $X$ une vari\'et\'e k\" ahl\'erienne
compacte, et $f:X\dasharrow X$ un endomorphisme m\'eromorphe
dominant. Il existe une application m\'eromorphe dominante $f$-invariante 
$g:X\dasharrow T$, dont la
fibre g\'en\'erale $X_t$ est l'adh\'erence de Zariski de l'ensemble
des it\'er\'es par $f$ d'un point g\'en\'eral de $X_t$.} 

\medskip

Ici comme dans le reste de l'article, le mot ``g\'en\'eral'' veut 
dire ``en dehors
d'une reunion d\'enombrable de ferm\'es analytiques stricts''.

\medskip

Les deux th\'eor\`emes pr\'ec\'edents impliquent facilement que, pour
$X={\cal F}(V)$ comme ci-dessus et assez g\'en\'eral, l'ensemble des
it\'er\'es par $f$ d'un point g\'en\'eral de $X$ est Zariski-dense
dans $X$.
En effet, $Pic(X)=\Z$; le th\'eor\`eme 2.1 garantit que 
$f^k, k>0$ ne pr\'eserve pas de fibration non-triviale (une
vari\'et\'e de type g\'en\'eral n'admet pas d'endomorphisme m\'eromorphe de
degr\'e $d>1$). Il s'ensuit, par la d\'ecomposition de Stein, que $dim(T)=0$
ou $dim(X)$ ($T$ comme dans le th\'eor\`eme 4.1); $dim(T)=dim(X)$ est 
impossible par $deg(f)>1$. Donc $T$ est un point.

\

Si l'on suppose $X$ projectif,  l'analogue du th\'eor\`eme 4.1 est 
vrai sur tout 
corps $K$ alg\'ebriquement clos de caract\'eristique
nulle. Si $K$ est d\'enombrable, cet enonc\'e n'a aucun contenu
non-trivial: il se peut qu'aucun
point de $X(K)$ n'est  ``g\'en\'eral'' au sens du th\'eor\`eme 4.1.
Par contre, nous pouvons appliquer le th\'eor\`eme \`a des
vari\'et\'es d\'efinies sur un corps non-d\'enombrable. 
Ainsi, un corollaire imm\'ediat de th\'eor\`emes 4.1 et
2.1
est la proposition suivante, motiv\'ee par \cite{HT2}:

\medskip

{\bf Proposition 4.8} {\it Soit $F$ un corps non-d\'enombrable (par
 exemple, le corps des fonctions m\'eromophes sur une courbe complexe).
 Soit $X$ une vari\'et\'e lisse projective
sur $\bar F$, v\'erifiant $K_X=0$, $Pic(X)=\Z$, et munie d'un
 endomorphisme 
rationnel
$f:X\dasharrow X$ de degr\'e $d>1$. Alors il existe une extension
finie
$L$ de $F$, telle que $X(L)$ soit Zariski-dense dans $X$.}

\medskip

La construction de C. Voisin ci-dessus fournit des exemples satisfaisant les
con- ditions de la proposition 4.8. 

\

Dans \cite{HT2}, des exemples analogues sont obtenus par une m\'ethode
diff\'erente. Notamment, les auteurs construisent des surfaces K3
de groupe de Picard cyclique, dont l'ensemble des points sur un
corps des fonctions est Zariski-dense. Comme nous l'avons d\'ej\`a
indiqu\'e, on ne sait pas si une telle surface admet des endomorphismes
rationnels de degr\'e $d>1$.

\

Puisque l'exemple de Voisin est hyperk\"ahl\'erien, nous sommes naturellement
amen\'es
\`a l'\'etude des fibrations
rationnelles (m\'eromorphes) sur 
les vari\'et\'es hy- perk\"ahl\'eriennes irr\'eductibles, lorsque la fibre g\'en\'erique est de dimension de
Kodaira nulle (section 3). Les fibrations r\'eguli\`eres ont \'et\'e
\'etudi\'ees
dans \cite{M1}, o\`u des r\'esultats tr\`es
pr\'ecis sont obtenus. 
Notre r\'esultat principal ici est le: 

\medskip 

{\bf Th\'eor\`eme 3.6} {\it Soit $X$ hyperk\"ahlerienne
  irr\'eductible projective
de dimension 4, et $f:X\dasharrow B$ une fibration 
\`a fibre g\'en\'erique de dimension
de Kodaira nulle. Alors $f$ poss\`ede un mod\`ele holomorphe, \`a un
flop pr\`es}.

\medskip 

Pr\'ecisons la phrase ``$f$ poss\`ede un mod\`ele holomorphe, \`a un
flop pr\`es'': elle veut dire qu'il existe une $X'$ hyperk\"ahlerienne
 irr\'eductible projective, et un flop
$\alpha: X\dasharrow X'$, telle que
$f\alpha^{-1}:X'\dasharrow B$ {\it poss\`ede un mod\`ele holomorphe},
c'est-\`a-dire qu'il existe une $\beta:B\dasharrow B^0$ birationnelle
avec $\beta f\alpha^{-1}$ holomorphe.

Apr\`es un flop, la situation est donc celle de
\cite{M1}. En particulier, les fibres de la fibration initiale sont
lagrangiennes, birationnelles \`a des surfaces ab\'eliennes. 

La preuve repose sur les 
r\'esultats de Wierzba (\cite{W}) et de Matsuki (\cite{Mat}), dont
on d\'eduit que, sur une $X$ comme ci-dessus, tout diviseur sans
composante 
fixe  peut \^etre
rendu nef \`a l'aide d'une suite finie de flops de Mukai.

\

Dans la section 5 nous avons rassembl\'e quelques remarques sur la
base d'une fibration m\'eromorphe
$g:X\dasharrow B$ avec $X$ hyperk\"ahlerienne irr\'eductible . 
Notamment, nous observons
que cette 
base $B$ n'a pas de tenseur holomorphe covariant non trivial. 
Cette propri\'et\'e entra\^{\i}ne la connexit\'e rationnelle 
de $B$ si $dim(B)\leq 3$, et conjecturalement en toute dimension.

\

Finalement, dans la section 6, nous situons l'\'etude des
endomorphismes $f:X\dasharrow X$ 
dans le cadre g\'en\'eral de la classification des vari\'et\'es k\" ahl\'eriennes compactes.

On montre (th\'eor\`eme 6.1) dans le cas g\'en\'eral ($X$
k\"ahlerienne compacte) 
qu'un it\'er\'e $f^N$ de $f$ pr\'eserve le {\it coeur} de $X$, 
introduit dans \cite{C3}. En particulier, si la $f$-orbite d'un
point est Zariski dense dans $X$, 
alors $X$ est sp\'eciale au sens de \cite{C3}. En g\'en\'eral, pour un certain
$N>0$, la fibration
$g_N:X\dasharrow T_N$, 
associ\'ee \`a $f^N$ par le th\'eor\`eme 4.1, factorise le coeur de $X$.

\

Quelques mois apr\`es avoir termin\'e la premi\`ere version  de ce texte
(math.AG/0510299), nous avons lu l'article \cite{M2} de D. Matsushita
(math.AG/0601114), dont le sujet est tr\`es li\'e au celui du n\^otre.
Ce travail de Matsushita nous a permis de g\'en\'eraliser
notre Proposition 3.1. Le th\'eor\`eme 3.6 peut, lui aussi, en
\^etre
d\'eduit; mais nous avons
n\'eanmoins pr\'ef\'er\'e exposer notre d\'emonstration initiale. En
effet, celle-ci est tr\`es \'el\'ementaire, et nous croyons que
notre point de vue peut avoir une certaine utilit\'e. 
Nous remercions D. Matsushita qui a attir\'e notre attention sur \cite{M2}.

\

Avec un tr\`es grand plaisir, nous d\'edions ce texte \`a Fedya Bogomolov,
auteur de travaux fondateurs sur les vari\'et\'es hyperk\"ahl\'eriennes et de nombreux autres sujets.
Ses articles en collaboration avec Yuri Tschinkel sur la densit\'e potentielle
de points rationnels ont aussi \'et\'e une source d'inspiration.

\

{\bf Notations et terminologie:}

\smallskip

On d\'enote par $\kappa(X)$ la dimension canonique, ou de Kodaira, de $X$, et par
$\kappa(X,L)$ la dimension d'Iitaka-Moishezon d'un fibr\'e en droites $L$.

Soient $X,B$ normales.
Une {\it fibration} est un morphisme surjectif $f: X\rightarrow B$ \`a fibres 
connexes. Une {\it fibration rationnelle} ou {\it m\'eromorphe} est
une application rationnelle dominante $g: X\dasharrow B$  \`a fibres 
connexes. Une {\it fibre} $X_b$ d'une fibration m\'eromorphe $g$
est  $\{\pi(x)| x\in \Gamma\ et\ f(x)=b\}$, o\`u $\Gamma$ est la
fermeture
du graphe de $g$ dans $X\times B$,
et $\pi:\Gamma \rightarrow X$, $f: \Gamma\rightarrow
B$ sont les deux projections. Parfois on appelle aussi une {\it fibre}
de $g$ une fibre de la fibration $Y\rightarrow B$, o\`u $Y$ est une
d\'esingularisation de $\Gamma$.

Nous dirons qu'une fibration est en vari\'et\'es de type
g\'en\'eral, en courbes, en sous-vari\'et\'es lagrangiennes, etc., 
si sa fibre g\'en\'erique
est une vari\'et\'e de type g\'en\'eral, une courbe, une
sous-vari\'et\'e
lagrangienne, etc.

Toutes les vari\'et\'es dans cet article sont compactes
k\"ahleriennes, et, dans les deux sections suivantes, elles sont, de
plus, suppos\'ees projectives.

Une vari\'et\'e {\it hyperk\"ahlerienne irr\'eductible}, ou  {\it symplectique
irr\'eductible}, est une vari\'et\'e lisse et simplement connexe $X$, telle que
$H^0(X,\Omega^2_X)$ est engendr\'e par une 2-forme holomorphe $\sigma$,
partout non-d\'eg\'en\'er\'ee sur $X$. Une telle vari\'et\'e est de
dimension paire $2n$, de classe canonique triviale, 
et si $n=1$, $X$ est une surface $K3$.

On notera parfois de la m\^eme fa\c con un diviseur de Cartier
et le fibr\'e en droites correspondant. Pour un tel diviseur $D$,
$|D|$ d\'enote le syst\`eme lin\'eaire complet associ\'e.

``G\'en\'eral'' veut dire ``en d\'ehors d'une r\'eunion d\'enombrable
de ferm\'es analytiques stricts; ``g\'en\'erique'' veut dire ``en d\'ehors
d'un ferm\'e analytique strict''.

Une fibration m\'eromorphe est dite presque holomorphe si sa fibre
g\'en\'erique ne rencontre pas le lieu d'ind\'etermination.

\section{ Fibrations sur les vari\'et\'es K-triviales g\'en\'erales}

Toute vari\'et\'e projective $X$ est birationnelle \`a
une fibration en diviseurs de base $\p^1$: il suffit de consid\'erer un pinceau de sections $P\subset
|L|$ d'un
fibr\'e
en droites $L$. Si $K_X=0$ et si $L$ est ample, les membres lisses de
$|L|$
sont de type g\'en\'eral. Cependant, il se peut que $P$ n'ait aucun
membre lisse.

\medskip

{\bf Th\'eor\`eme 2.1} {\it Soit $X$ une vari\'et\'e lisse projective,
telle que $K_X= 0$ et $NS(X)\cong \Z$. Alors toute fibration rationnelle
$g:X \dasharrow B$, $0<dim(B)<dim(X)$, est en vari\'et\'es de type 
g\'en\'eral.}

\medskip

{\it Preuve:} Soit $g$ une telle fibration. La r\'esolution des
singularit\'es fournit
alors une
vari\'et\'e lisse projective $Y$, une fibration $f:Y\rightarrow B$
et
un morphisme birationnel $\pi: Y\rightarrow X$:
$$
\renewcommand{\arraystretch}{0}
\arraycolsep=0.05em
\begin{array}{ccccc}
&& Y\\
& \stackrel{\pi}{\swarrow} && \stackrel{f}{\searrow} &\\
X && \stackrel{g}{\dasharrow} && B
\end{array}
$$

Soit $H$ un diviseur tr\`es ample sur $B$ et $L=\pi_*f^*H$ (on
 rappelle
qu'ici,  $\pi_*$
est consid\'er\'e au sens de cycles; comme $X$ est lisse, $L$ est un diviseur de
Cartier), de sorte que l'application $g$ est donn\'ee par un
 syst\`eme
lin\'eaire $U\subset |L|$. On a alors
$$\pi^*L=f^*H + \sum_ia_iE_i,$$
o\`u les $E_i$ sont les diviseurs $\pi$-exceptionnels et $a_i\geq
0$.
D'autre part, puisque $X$ est lisse et K-triviale,

$$K_Y=\sum_ie_iE_i$$
avec $e_i$ strictement positifs. La restriction \`a la fibre
g\'en\'erique $F$ de $f$ donne 
$$\pi^*L|_F=\sum_i a_iE_i|_F;\ K_F=\sum_ie_iE_i|_F,$$
puisque le fibr\'e normal de $F$ dans $Y$ est trivial.
Mais, comme tous les $e_i$ sont strictement positifs, il existe
un entier positif $m$ tel que $\pi^*L|_F<mK_F$, c'est-\`a-dire,
$mK_F-\pi^*L|_F$ est effectif. Ceci implique 
$$\kappa(F, \pi^*L|_F)\leq \kappa(F).$$

Rappellons maintenant que $NS(X)=\Z$, et que, par construction, 
$L$ est un diviseur
effectif et
non-nul. Donc $L$ est ample, 
$\kappa(F, \pi^*L|_F)=dim(F)$, et $\kappa(F)=dim(F)$.

\medskip

{\bf Corollaire 2.2} {\it  Soit $X$ comme dans le th\'eor\`eme 2.1, et
  $h:X\dasharrow X$
un endomorphisme rationnel de degr\'e $>1$. Alors $h$ ne pr\'eserve
  pas de
fibration non-triviale.}

\medskip

En effet, les vari\'et\'es de type g\'en\'eral ne poss\`edent pas
d'endomorphismes rationnels de degr\'e $>1$.

\

{\bf Variantes et remarques:}

\medskip

1) L'argument de la d\'emonstration du th\'eor\`eme 2.1 montre aussi
   la variante suivante:

\smallskip

{\it Soit $X$ telle que $\kappa(X)\geq 0$, $L$ un fibr\'e en
  droites
sur $X$, et $U\subset |L|$ un syst\`eme lin\'eaire d\'efinissant une
fibration rationnelle $g: X\dasharrow B$. Si $F$ est une fibre
g\'en\'erique de $g$, alors}

$$\kappa(X, L)\leq dim(B)+\kappa(F)\leq dim(U)+\kappa(F).$$

De plus, un l\'eger raffinement de cet argument permet de pr\'eciser le th\'eor\`eme 2.1 
comme suit:

\medskip

{\bf Th\'eor\`eme 2.3} {\it  Soit $X$ compacte k\"ahl\'erienne telle que $K_X=0$.
Soit $g:X\dasharrow B$ une fibration m\'eromorphe avec $B$ projective, et $L:=\pi_
*(f^*(H))\in Pic(X)$, $H$ ample sur $B$. 
Si $F$ est la fibre g\'en\'erique d'une r\'esolution $f:Y\to B$, alors 
$\kappa(X,L)=dim(B)+\kappa(F)$. }

\medskip

{\it D\'emonstration:} On conserve les notations introduites dans la
d\'emonstration de 2.1. Soit $E'$ un $\Q$-diviseur effectif et $\pi$-exceptionnel
arbitraire sur $Y$. Alors $\kappa(Y,\pi^*(L)+E')=\kappa(X,L)=\kappa(Y,\pi^*(L))$,
par le th\'eor\`eme d'Hartogs (qui permet d'\'etendre \`a $X$ les sections de $mL$, $m>0$ entier, 
d\'efinies sur l'ouvert au-dessus duquel $\pi$ est un isomorphisme). On peut choisir $E'$ tel que 
$\pi^*(L)+E'=f^*(H)+D$ o\`u $D=\sum_i d_iE_i$ avec $d_i>0$ pour tout $i$, de sorte que $D$
contient tout diviseur $\pi$-exceptionnel. Comme dans 2.1, on a $\kappa(F)=\kappa(F, D|_F)
=\kappa(F,(f^*(H)+D)|_F)$. La conclusion r\'esulte donc du lemme suivant: 

\medskip

{\bf Lemme 2.4} {\it Soit $f:Y\to B$ une fibration holomorphe, avec 
$Y$ une vari\'et\'e complexe compacte et $B$ projective. 
Soit $H$ un fibr\'e en droites ample sur $B$ et $D$ un $\Q$-diviseur effectif sur $Y$. 
Alors, pour $F$ la fibre g\'en\'erique de $f$,  $$\kappa(Y,f^*(H)+D)=dim(B)+\kappa(F, D|_F)$$.}

{\it Preuve:} Puisque $D$ est effectif, la r\'eduction d'Iitaka $I: Y\rightarrow Z$ de $f^*(H)+D$ factorise
$f$. On en d\'eduit que la fibre g\'en\'erique de $I$ coincide avec celle de $I_F$, la r\'eduction d'Iitaka
de $(f^*(H)+D)|_F=D|_F$. Lemme 2.5 en r\'esulte imm\'ediatement.

\medskip

2) Ce m\^eme argument s'applique aussi lorsque $X$
est $\Q$-factorielle \`a singularit\'es terminales (le diviseur $L$ 
est alors $\Q$-Cartier, et son image inverse et sa dimension
d'Iitaka-Moishezon
 sont donc d\'efinies; et les ``discr\'epances'' des diviseurs
exceptionnels
sont toujours strictement positives). 

Cependant, il existe des
vari\'et\'es $K$-triviales \`a points doubles, \`a groupe de
Picard cyclique, admettant une fibration rationnelle non triviale en
vari\'et\'es
de $\kappa$ nulle: telle est la {\it quintique de
Horrocks-Mumford} (\cite{Bor}). Nous remercions
B. Hassett pour nous avoir signal\'e cet exemple.

\smallskip

3) Une vari\'et\'e $K$-triviale peut \^etre non trivialement
domin\'ee 
par une famille de
   vari\'et\'es de dimension canonique nulle: par exemple, toute
surface K3 est (g\'en\'eriquement) recouverte par des courbes
   elliptiques 
(en g\'en\'eral singuli\`eres, \cite{MM}). Si le  groupe de Picard est cyclique, une telle 
famille ne fournit cependant jamais de fibration m\'eromorphe.

\section{ Fibrations sur les  vari\'et\'es hyperk\"ahl\'eriennes}

Les exemples de C. Voisin (mentionn\'es dans l'introduction)
de dimension minimale (quatre) sont hyperk\"ahl\'eriens irr\'eductibles. 
Il est donc int\'eressant de
comprendre
les fibrations m\'eromorphes sur les vari\'et\'es 
hyperk\"ahleriennes irr\'eductibles, dont les fibres ne sont pas de type
g\'en\'eral. A l'aide de la fibration d'Iitaka-Moishezon
r\'elative, on se ram\`ene aux fibrations m\'eromorphes dont les
fibres g\'en\'eriques sont de dimension canonique nulle. 

De tels exemples existent pour
les
vari\'et\'es plus sp\'eciales (de nombre de Picard $2$ au moins), tels $Hilb^n(S)$ pour certaines
surfaces
$S$ de type K3 (voir par exemple \cite{HT}). 
Les exemples connus sont des fibrations en
tores lagrangiens de dimension moiti\'e.

Dans \cite{M1}, il est d\'emontr\'e que si $X$ est hyperk\"ahl\'erienne 
irr\'eductible projective de dimension $2n$ et 
$f: X\rightarrow B$ est une fibration r\'eguli\`ere, alors la fibre g\'en\'erique
de $f$ est un tore lagrangien de dimension $n$, $Pic(B)=\Z$ et $B$ est
$\Q$-Fano \`a singularit\'es log-terminales $\Q$-factorielles.

\

On peut esp\'erer que si $g: X\dasharrow B$ est seulement rationnelle 
\`a fibres de
dimension
de Kodaira nulle, la situation reste
semblable \`a celle de \cite{M1}. 

\

{\bf Question 3.0:} {\it Si $g: X\dasharrow B$ est une fibration rationnelle non 
triviale, avec $X$ hyperk\" ahl\'erienne irr\'eductible et projective, \`a 
fibre g\'en\'erique $F$ telle que $\kappa(F)=0$, alors les propri\'et\'es 
suivantes sont-elles satisfaites?

1. $2dim(B)=dim(X)$

2. $F$ est lagrangienne, et birationnelle \`a une vari\'et\'e ab\'elienne.}

\

Montrons d'abord qu'une telle fibration est donn\'ee par des sections 
d'un fibr\'e en droites $L$ tel que $L^{dim(X)}=0$.

\

Soit $g: X\dasharrow B$ comme ci-dessus. On fixe $H\in Pic(B)$ engendr\'e
par ses sections, et soit
$L=g^*H$. Plus pr\'ecisement, comme dans la premi\`ere section,
$L=\pi_*f^*H$, o\`u
$$
\renewcommand{\arraystretch}{0}
\arraycolsep=0.05em
\begin{array}{ccccc}
&& Y\\
& \stackrel{\pi}{\swarrow} && \stackrel{f}{\searrow} &\\
X && \stackrel{g}{\dasharrow} && B
\end{array}
$$
est une r\'esolution du lieu d'ind\'etermination de $g$ (avec $Y$
lisse). 
On a alors
$\pi^*L=f^*H+\sum a_iE_i$ avec $E_i$ $\pi$-exceptionnels, $a_i\geq 0$.

\

On d\'enote par $\sigma$ la forme symplectique sur $X$.
Soit $q$ la forme quadratique de Beauville-Bogomolov sur $H^2(X, \R)$.
On rappelle (voir par exemple \cite{H1}) que pour $\alpha$ de type
$(1,1)$,
$$q(\alpha)=c\int \alpha^2(\sigma\overline{\sigma})^{n-1},$$
o\`u $c$ est une constante positive, et qu'il existe une autre
constante $c'$ telle que pour tout $\alpha$,
$\alpha^{2n}=c' q(A)^n$.

\medskip

{\bf Proposition 3.1} {\it 1) La classe du fibr\'e en droites $L$ dans 
$H^2(X, \R)$ est isotrope par rapport \`a la forme $q$. Par cons\'equent,
$L^{dim(X)}=0$.

2) Soit maintenant $g:X\dasharrow B$ une application donn\'ee par un
syst\`eme lin\'eaire de sections d'un fibr\'e en droites $L$, $q(L)=0$.
Alors la fibre g\'en\'erique de $g$ n'est pas de type g\'en\'eral.}

\medskip

{\it Preuve:} Remarquons d'abord que $L$ est sans composantes fixes,
d'o\`u $q(L)\geq 0$. En effet, soient $L_1, L_2$ deux membres irr\'eductibles 
distincts de
$\vert L\vert$;
$V=L_1\cap L_2$ est donc de codimension pure deux dans $X$. On a 
$$q(L)=c\int L^2(\sigma\overline{\sigma})^{n-1}=
c\int_V(\sigma\overline{\sigma})^{n-1}\geq 0$$ puisque 
$(\sigma\overline{\sigma})^{n-1}$ est positive.

Supposons $q(L)>0$. Alors, selon \cite{Bo}, th\'eor\`eme 4.3 (i), $L$ est
dans l'int\'erieur du c\^one pseudo-effectif de $X$, donc c'est un fibr\'e
``vaste'', ou ``grand'' (traduction de {\it big}). Il s'ensuit que pour la fibre g\'en\'erique $F$ de $f$,
$\pi^*L|_F$ est vaste aussi. Comme on l'a d\'ej\`a remarqu\'e dans la
d\'emonstration du th\'eor\`eme 2.1, pour un certain entier positif $m$,
$mK_F-\pi^*L|_F$ est effectif. Donc $\kappa(F)\geq
\kappa(F,\pi^*L|_F)=dim(F)$,
 autrement dit,
$F$ est de type g\'en\'eral, une contradiction.

R\'eciproquement, \'etant donn\'e un $L$ sans composantes fixes, tel 
que $q(L)=0$,
on voit que $\kappa(L)<dim(X)$ (nous remercions S. Boucksom
pour cette remarque): en effet, 
$q(L,M)\geq 0$ pour tout $M$ effectif. Si $L$ est vaste, alors $L=A+B$
avec $A$ ample et $B$ effectif. Par consequent, $q(L)\geq q(A) >0$.
Le second enonc\'e du proposition 3.1 r\'esulte donc de 2.3.

\medskip

{\it Remarque:} Cet argument est inspir\'e par \cite{M2}. Notre
argument initial marchait modulo le programme de mod\`eles minimaux (MMP) 
pour les 
vari\'et\'es de
dimension de Kodaira nulle en dimension $\leq dim(X)-1$. Le MMP servait \`a
trouver, par un point g\'en\'erique de $X$, une courbe $C$ telle que $LC=0$.
On en d\'eduisait $q(L)=0$ \`a l'aide de \cite{H2}.

\

Le corollaire suivant est motiv\'e par les exemples de Hassett et
Tschinkel dans \cite{HT}; ils d\'emontrent que  $Hilb^2(S)$, o\`u $S$ est 
une surface
$K3$
g\'en\'erique de degr\'e $2m^2, m>1$, est isomorphe \`a une fibration
en
tores au-dessus de $\p^2$ (si $m=1$, il est clair que 
$Hilb^2(S)$ est birationnelle \`a une telle
fibration; en effet, $S$ est un rev\^etement double $h: S\rightarrow
\p^2$, et $h$ induit $g: Hilb^2(S)\dasharrow (\p^2)^*$, o\`u pour $Z$ 
g\'en\'erique, $g(Z)$ est la droite de $\p^2$ contenant $h(Z)$).

\medskip

{\bf Corollaire 3.2} {\it Soit $S$ une surface K3 g\'en\'erique de
  degr\'e
$d\neq 2m^2, m\in \Z$. Alors $Hilb^2(S)$ n'est pas birationnelle \`a une
fibration en vari\'et\'es de dimension canonique nulle.}

\medskip

En effet, si $d\neq 2m^2$, alors la forme de 
Beauville-Bogomolov ne repr\'esente
pas z\'ero sur $Pic(Hilb^2(S)).$

\

Par cons\'equent, toute fibration m\'eromorphe non-triviale de
$Hilb^2(S)$ comme dans le corollaire, est en vari\'et\'es de type 
g\'en\'eral.

\

La proposition suivante est l'analogue du fait que, dans le cas o\`u
$g:X\to B$ est holomorphe, $Pic(B)=\Z$ (\cite{M1}): 

\medskip

{\bf Proposition 3.3} {\it Pour tout $M\in Pic(B)$,
$\pi_*f^*M$ est proportionnel \`a $L$.}

\medskip

{\it Preuve:} Il suffit de d\'emontrer l'enonc\'e pour $M$ tr\`es ample, 
puisque de tels $M$ engendrent $Pic(B)$. Soit $N=\pi_*f^*M$. 
On montre, comme
dans la proposition pr\'ec\'edente, que $q(N)=0$; de plus, pour 
$a,b$ assez grands, $|aH+bM|$ est tr\`es ample, d'o\`u
$q(aL+bN)=0$ pour $a,b$ assez grands, et par cons\'equent pour
$a,b$ arbitraires.

On sait (\cite{M1}) que, pour $A\in Pic(X)$ ample et $B\in Pic(X)$
$q$-isotrope,
$BA^{2n-1}=0$ implique $B=0$. Fixons $A\in Pic(X)$ ample; alors,
pour un certain $\lambda\in \Q$, $(\lambda L-M)A^{2n-1}=0$, d'o\`u
la
proposition.

\

Voici encore une application imm\'ediate de la proposition 3.1:

\medskip

{\bf Proposition 3.4} {\it La dimension de la fibre g\'en\'erique $F$ de $g$ 
est $\geq n$ (o\`u $2n=dim(X)\geq 4$). }

\medskip

{\it Preuve:} Supposons le contraire.
Soit $b=f(F)$, et soit
$V\subset T_bB$ un sous-espace de codimension 2. Rempla\c{c}ons, si
n\'ecessaire, $H$ par $lH$ avec $l$ assez grand, et
consid\'erons deux diviseurs
$H_1, H_2\in |H|$ passant par $b$, tels que $V=T_b(H_1\cap H_2)$.
Soient $L_i$ les transform\'ees strictes des $H_i$ sur $X$;
comme $0=q(L)=c\int_{L_1\cap L_2}(\sigma\overline{\sigma})^{n-1}$,
la restriction de $\sigma^{n-1}$ \`a toute composante irr\'eductible de $L_1\cap L_2$
est nulle (sinon l'int\'egrale serait strictement positive). Donc la restriction de $\sigma$
\`a une telle composante est d\'eg\'en\'er\'ee.
Soit $c\in C=\pi(F)$ g\'en\'erique; donc $C$ est lisse en $c$ et $\pi$ est un isomorphisme
au voisinage de $c$. On en d\'eduit que pour tout sous-espace $W$ de codimension 2 
dans
$T_cX$, contenant $T_cC$, la restriction de $\sigma$ \`a $W$ est d\'eg\'en\'er\'ee.
Un calcul facile d'alg\`ebre lin\'eaire montre alors que c'est
impossible, puisque $\sigma$ est non-d\'eg\'en\'er\'ee.

\

Nous allons ensuite \'etudier en d\'etail le cas $dim(X)=4$.

On se place donc, jusqu'\`a la fin de cette partie, dans la situation suivante:

\

{\bf ($\ast$)} {\it Soit $g: X\dasharrow B$ une fibration rationnelle non
triviale, avec $X$ hyperk\" ahl\'erienne irr\'eductible et
projective de dimension $4$, \`a fibre g\'en\'erique $F$
telle que $\kappa(F)=0$. On garde les notations d\'ej\`a introduites
($Y$, $f$, $\pi$ etc.). $H$ d\'esigne un diviseur tr\`es ample sur $B$,
et $L=\pi_*f^*H$.} 

\

On fera appel au r\'esultat suivant de
Wierzba (\cite{W}):

\

{\bf Th\'eor\`eme (J. Wierzba)} {\it Soit $X$ une vari\'et\'e
hyperk\"ahlerienne irr\'eductible projective de dimension 4 et $D$ 
un diviseur sans
composante
fixe sur $X$. Il existe une autre vari\'et\'e
hyperk\"ahlerienne irr\'eductible projective $X'$ de dimension 4
et une application birationnelle $\phi: X\dasharrow X'$, telle
que
la transform\'ee stricte $D'=\phi_*D$ est nef sur $X'$. De plus,
$\phi$ est un produit fini de flops de Mukai (projectifs).}

\

{\it Remarque:} B. Fu nous a signal\'e que la preuve de \cite{W} 
est incompl\`ete;
mais que l'erreur se trouve dans l'argument de la terminaison de flops.
Puisque cette terminaison est d\'emontr\'e en dimension 4 par
Matsuki (\cite{Mat}), nous pouvons utiliser ce r\'esultat.

\

Soit donc $g:X\dasharrow B$ notre fibration donn\'ee par un
syst\`eme lin\'eaire $U\subset |L|$; consid\'erons $\phi: X\dasharrow
X'$
comme dans le th\'eor\`eme de Wierzba; $\phi$ est un isomorphisme
en codimension un. La
transform\'ee
stricte $L'=\phi_*L$ est nef; on a un syst\`eme lin\'eaire
$U'=\phi_*U\subset |L'|$ qui d\'etermine une fibration rationnelle
$g':X'\dasharrow B'$. La fibre g\'en\'erique de $g'$ est
birationnelle
\`a celle de $g$. De plus, $L'$ est isotrope pour la
forme
de Beauville-Bogomolov.

\

Autrement dit, pour comprendre la fibre g\'en\'erique $\pi(F)$ de $g$, 
on peut supposer que $L$
est nef.

\medskip

{\bf Lemme 3.5} {\it Dans la situation {\bf $(\ast)$}, 
$g$ est presque holomorphe si $L=\pi_*f^*H$ est nef.}

\medskip

{\it Preuve:} Ecrivons $\pi^*L=f^*H+\sum_ia_iE_i,\ a_i\geq 0$; 
$\pi^*L$ est nef,
donc $\sum_ia_iE_i$ est $f$-nef: $(\sum_ia_iE_i|_F)\cdot C\geq 0$
pour tout courbe $C$ dans une fibre de $f$. En particulier, $\sum_ia_iE_i|_F$
est nef sur $F$. Mais cette somme
est support\'ee par le diviseur canonique de $F$, et $\kappa(F)=0$.
Comme $dim(X)=4$, on a  $dim(F)=2$ ou $3$; $F$ poss\`ede donc
un mod\`ele minimal $F_0$, et il y a un diagramme

$$
\renewcommand{\arraystretch}{0}
\arraycolsep=0.05em
\begin{array}{ccccc}
&& \tilde{F}\\
& \stackrel{\mu}{\swarrow} && \stackrel{\tilde{h}}{\searrow} &\\
F && \stackrel{h}{\dasharrow} && F_0
\end{array}
$$
avec $\mu$ et $\tilde{h}$ des morphismes birationnels.
Le diviseur $\mu^*(\sum_ia_iE_i|_F)$ est nef, support\'e par une
partie
du diviseur canonique de $\tilde{F}$. Le diviseur canonique de
$\tilde{F}$ est contract\'e par $\tilde{h}$. Un diviseur
nef et contractible est nul. Donc
$\mu^*(\sum_ia_iE_i|_F)$ est nul, et $E_i\cap
F=\emptyset$ d\`es que $E_i$ intervient dans $\pi^*L$. Donc
$g$ est presque holomorphe: en effet, aucune courbe $C$
$\pi$-exceptionnelle
mais non $f$-exceptionnelle ne rencontre $F$ (puisque
$(\sum_ia_iE_i)C<0$),
ce qui implique que $\pi(F)$ ne rencontre pas le lieu d'ind\'etermination.

\

{\bf Remarques:} 1) L'implication ``si $L$ est nef, alors $g$ est
presque holomorphe'' est vraie en toute dimension modulo le 
programme de mod\`eles minimaux (MMP). De plus, dans ce cas, MMP ram\`ene
notre situation a celle de \cite{M2}: $\kappa(F)=0$ implique alors
que
la dimension num\'erique de $L$ n'est pas maximale. Comme on l'a
d\'ej\`a dit dans l'introduction, nous allons cependant continuer
notre d\'emonstration, qui peut pr\'esenter un int\'er\^et
ind\'ependant.

2) On voit d\'ej\`a facilement que la fibre $F$ est de dimension 2
(la proposition 3.4 et la presque-holomorphie de $F$ emp\^echant 
$dim(B)=1$) et que $F$ est
lagrangien (l'argument de 3.4), et est par cons\'equent un tore (th\'eor\`eme
de Liouville sur les syst\`emes int\'egrables, \cite{Arn}, p. 271).
Nous allons d\'emontrer un r\'esultat plus pr\'ecis:

\

{\bf Th\'eor\`eme 3.6} {\it Si $L$ est nef, $g$ poss\`ede un mod\`ele
holomorphe, c'est-\`a-dire qu'il existe une application birationnelle
$\psi: B\dasharrow B_0$, telle que $\psi g$ est r\'eguli\`ere. Ainsi,
toute fibration $X\dasharrow B$ \`a fibre g\'en\'erique de dimension
de Kodaira nulle poss\`ede un mod\`ele
holomorphe \`a un flop pr\`es}.

\medskip

Par cons\'equent, la r\'eponse \`a la question 3.0 est positive en
dimension 4: le cas holomorphe
\'etant vrai par \cite{M1}.

\

Avant de commencer la d\'emonstration, \'etablissons un lemme
technique facile:

\medskip

{\bf Lemme 3.7} {\it Soit $f:X\rightarrow B$ une fibration
  holomorphe (on suppose 
$X$, $B$ normales projectives mais pas n\'ec\'essairement lisses; a fortiori,
on ne suppose pas, dans ce lemme, que $X$ est symplectique de dimension 4). 
Soit $E$ un diviseur
  de
Cartier effectif sur $X$, ne dominant pas $B$ et $f$-nef. Alors, 
pour chaque composante irr\'eductible
$A\subset f(Supp(E))$ de codimension 1 dans $B$, $Supp(E)$ 
contient $f^{-1}(a)$ pour $a\in A$ 
g\'en\'erique.}

\medskip

{\it Preuve:} Soit $r: X'\rightarrow X$ une
d\'esingularisation. Alors
$r^*(E)$ satisfait \`a la condition du lemme, et la conclusion pour 
$r^*(E)$ entraine celle pour $E$; on peut donc supposer que $X$ est
lisse. $B$ \'etant normale, une courbe g\'en\'erique ample $H$
 sur $B$ est lisse, ainsi que $V=f^{-1}H$, qui
est une fibration au-dessus de $H$ dont toutes les fibres ont la m\^eme
dimension.
Alors $H$ coupe toute composante $A\subset Supp(f(E))$ de codimension 1 dans $B$ 
au point g\'en\'erique
de $A$.
Soit $S$ une surface ample assez g\'en\'erale sur $V$;
$f|_S:S\rightarrow H$ est une fibration d'une surface lisse sur une
courbe
lisse. La restriction de $E$ sur $S$ est nef, et support\'ee sur des
fibres
de $f|_S$. Le lemme de Zariski (\cite{BPV}, p.90) implique que
$E|_S=f|_S^*D$ avec $D\in Pic(H)$. Donc $E|_S$ contient toute
composante
irr\'eductible de toute fibre de $f_S$ qui le supporte; par
construction,
la m\^eme chose est vraie pour $E$ au-dessus d'un point g\'en\'erique de $A$. 

\

{\it D\'emonstration du th\'eor\`eme 3.6:} 

Ecrivons $\pi^*L=f^*H+E,$

$$
\renewcommand{\arraystretch}{0}
\arraycolsep=0.05em
\begin{array}{ccccc}
&& Y\\
& \stackrel{\pi}{\swarrow} && \stackrel{f}{\searrow} &\\
X && \stackrel{g}{\dasharrow} && B,
\end{array}
$$
o\`u $E$ est $\pi$-exceptionnel effectif et $H$ est tr\`es ample.  
Le diviseur $E$ a des composantes de deux
types
diff\'erents: celles qui se projettent sur une courbe de $B$, et celles
dont l'image par $f$ est un point
(autrement dit, support\'ees sur des fibres de dimension 3 de $f$).

Ecrivons $E=E'+E''$ de sorte que l'ensemble $f(E'')$ est fini, 
et $f$ est de dimension relative pure \'egale \`a 2
sur
$Supp(E')$. On a la d\'ecomposition en composantes irr\'eductibles 
$Supp(E')=\cup E'_i$; soit $D_i = f(E'_i)$. 

Nous aimerions construire
une application birationnelle de $B$ dans la vari\'et\'e de Chow de
$X$ qui contracte les $D_i$, en faisant correspondre le cycle 
$\pi_*([Y_b])$
\`a  $b\in B$. Nous avons un obstacle pour ce faire: les fibres de
$f$ ne sont pas forc\'ement toutes de dimension 2.
 Apr\`es applatissement g\'eom\'etrique et normalisation, on a 
 un diagramme commutatif

$$
\renewcommand{\arraystretch}{0}
\arraycolsep=0.05em
\begin{array}{cccccccccc}
&& Y && \stackrel{\tilde{\pi}}{\longleftarrow} && \tilde{Y} \\
& \stackrel{\pi}{\swarrow} && \stackrel{f}{\searrow} &&&& \stackrel{\tilde{f}}{\searrow}\\
X && \stackrel{g}{\dasharrow} && B && \stackrel{\beta}{\longleftarrow} && \tilde{B},
\end{array}
$$

avec $\beta$ birationnel, $\tilde{B}$ normal, $\tilde{Y}$ la
normalisation de la composante irr\'eductible de $Y\times_B \tilde{B}$ qui domine
$\tilde{B}$ (en g\'en\'eral, il y a d'autres composantes, ne dominant pas $\tilde{B}$,
qui proviennent de
fibres
de dimension 3 quand on \'eclate la base au point correspondant), et 
toutes les
fibres de $\tilde{f}$ de dimension 2. Par \cite{Bar2}, nous avons un morphisme 
$p:\tilde{B}\rightarrow Chow(X)$,
$\tilde{b}\mapsto (\pi\tilde{\pi})_*([\tilde{Y}_{\tilde{b}}])$. 

\ 

Soit $B_0=Im(p)$; nous affirmons que l'application
$p{\beta}^{-1}g: X\dasharrow B_0$ est holomorphe. 

Regardons d'abord l'application ${\beta}^{-1}g: X\dasharrow \tilde{B}$.
L'\'egalit\'e $\pi^*L= f^*H+E$ fournit 
$$\tilde{\pi}^*\pi^*L = \tilde{\pi}^* f^*H+\tilde{\pi}^*E = 
(\tilde{f})^*\beta^*H+\tilde{\pi}^*E,$$ de sorte que 
$L=(\pi\tilde{\pi})_*(\tilde{f})^*(\beta^*H).$ Le diviseur
$\beta^*H$ n'est pas, en g\'en\'eral, tr\`es ample, on ne peut donc
pas affirmer que ${\beta}^{-1}g$ est donn\'ee par des sections de
$L$;
mais, d'apr\`es la proposition 3.3, pour un $\tilde{H}$ tr\`es ample sur
$\tilde{B}$,
$M=(\pi\tilde{\pi})_*(\tilde{f})^*\tilde{H}$ est proportionnel \`a
$L$,
donc \`a nouveau nef. Nous avons
$(\pi\tilde{\pi})^*M=(\tilde{f})^*\tilde{H}+\tilde{E}$ o\`u
$\tilde{E}$ est $\pi\tilde{\pi}$-exceptionnel et ne domine pas
$\tilde{B}$.
De plus, il est $\tilde{f}$-nef parce que $M$ est nef. Toute
composante $\tilde{E}_i$
de $\tilde{E}$ se projette par $\tilde{f}$ sur une courbe $\tilde{D}_i\subset
\tilde{B}$. 

On applique maintenant le lemme 3.7: du moins au point g\'en\'erique
$d$ de chaque
$\tilde{D}_i$,
la fibre de $\tilde{f}$ est toute enti\`ere dans le support de
$\tilde{E}$.

Ceci veut dire que $p(d)$ est un cycle support\'e par le lieu
exceptionnel $Z$
de $\pi\tilde{\pi}$. Comme $Z$ est de dimension deux, l'ensemble des
valeurs possibles de $p(d)$ est discret. Donc $p(d)$
est
constant quand $d$ parcourt $\tilde{D}_i$; autrement dit, $p$ contracte les
$\tilde{D}_i$.

Soit maintenant $x$ un point d'ind\'etermination de
$p{\beta}^{-1}g$. Alors
il existe une courbe irr\'eductible $C$ sur $\tilde{Y}$, telle que
$\pi\tilde{\pi}(C)=x$ et $C$ n'est pas contract\'ee par $p\tilde{f}$.
En particulier, $C$ n'est pas contract\'ee par $\tilde{f}$.
Comme $(\pi\tilde{\pi})^*M\cdot C = 0$ et $\tilde{H}$ est tr\`es
ample,
on doit avoir $\tilde{E}\cdot C < 0$, donc $C\subset
Supp(\tilde{E})$.
Par consequent, $\tilde{f}(C)=\tilde{D}_i$ pour un certain $i$, et donc $f(C)$ est
contract\'e par $p$, ce qui est absurde. Nous avons donc d\'emontr\'e
que le lieu d'ind\'etermination de $p{\beta}^{-1}g$ est vide,
autrement dit,  $p{\beta}^{-1}g$ est holomorphe. On peut donc
prendre $p{\beta}^{-1}$  pour $\psi$,  C. Q. F. D.

\

{\bf Remarques:} 

1) Une version d'une conjecture bien connue depuis quelques ann\'ees dit
qu'une fibr\'e en droites $L$ nef, non-nul et tel que $q(L)=0$ d\'efinit
une fibration lagrangienne (holomorphe) d'une $X$ hyperk\"ahlerienne. Ici,
nous l'avons d\'emontr\'e en dimension 4, avec une hypoth\`ese suppl\'ementaire
que $L$ a du moins un faisceau de sections.

2) L'argument
sur la vari\'et\'e de Chow
se g\'en\'eralise-t-il en dimension quelconque? Il suffirait, pour une telle 
g\'en\'eralisation, d'avoir $dim(I(g))\leq dim(X)/2$, avec $I(g)$ le lieu
d'ind\'etermination de $g$, lorsque $L$ est nef. Une telle
g\'en\'eralisation
\'eviterait d'invoquer le r\'esultat
difficile de \cite{Kaw} dans \cite{M2}.

3) Si $g:X\to B$ est holomorphe, $B$ est une surface
$\Q$-Fano \`a nombre de Picard $1$, \`a singularit\'es
log-terminales et $\Q$-factorielles ([M1]). On peut se demander s'il
s'agit toujours de $\p^2$. On peut v\'erifier facilement que c'est le
cas si toute fibre de $g$ a une composante r\'eduite. 
De plus, le lieu lisse de $B$ est simplement connexe. Effectivement,
puisque toutes les fibres de $g$ sont de dimension 2 par \cite{M3}, 
$g^{-1}(Sing(B))$ est de codimension 2 dans $X$ (si non-vide), donc 
$X-g^{-1}(Sing(B))$
est, tout comme $X$, simplement connexe; ce qui, par la connexit\'e des fibres
de $g$, implique que $B-Sing(B)$ est simplement connexe.

 \section{Fibration associ\'ee \`a un endomorphisme} 

Soit $X$ une vari\'et\'e k\"ahl\'erienne compacte et connexe
de dimension complexe $n$, et $f:X\dasharrow X$ un endomorphisme m\'eromorphe
dominant de $X$. Pour $x\in X$ g\'en\'eral (i.e: dans une intersection
d\'enombrable d'ouverts de Zariski denses de $X$), l'application $f^m$ est
holomorphe en $x$ pour tout entier $m>0$. On appelle alors l'ensemble des
$f^m(x), m\geq 0$ la $f$-orbite de $x$, et on note $Z^+_f(x)$
(ou
simplement $Z^+(x)$)
l'adh\'erence de Zariski de cette $f$-orbite dans $X$. (Par un ferm\'e (resp. un ouvert) de 
Zariski sur $X$, nous entendons un ferm\'e (resp. un ouvert) analytique. Lorsque $X$ est 
projective, c'est la topologie de Zariski ``usuelle''.)
 
\medskip

On se propose de d\'emontrer le r\'esultat
suivant:

\medskip

{\bf Th\'eor\`eme 4.1} 
{\it Dans la situation pr\'ec\'edente, il existe une application m\'eromorphe 
dominante
 $g:X\dasharrow T$ telle que:
 
1. $g\circ f=g$

2. Pour $x\in X$ g\'en\'eral, la fibre de $g$ passant par $x$ est \'egale
\`a $Z^+_f(x)$.}

\medskip

En appliquant la d\'ecomposition de Stein, on obtient une fibration
m\'eromorphe (possiblement triviale)  pr\'eserv\'ee par une certaine puissance de $f$.

\

{\bf {Remarques:}}

1. Si en un point $a$ de $X$, toutes les it\'er\'ees $f^m(a)$ de $a$ par $f$ sont d\'efinies, et si elles sont Zariski-denses dans $X$, il en est donc de m\^eme pour un point g\'en\'eral de $X$.

2. Si $X$ et $f$ sont d\'efinis sur un corps de nombres, et si la $f$-orbite d'un point $a$ de $X$ (au sens pr\'ec\'edent) est Zariski-dense, on aimerait savoir si l'on peut choisir $a$ d\'efini sur un corps de nombres (auquel cas $X$ serait potentiellement dense). Ceci semble plausible, mais ne peut pas \^etre d\'emontr\'e par les m\'ethodes utilis\'ees ici, puisque les points de $X$ d\'efinis sur $\bar \Q$ ne sont pas g\'en\'eraux. 

3. Si $f$ est de degr\'e $d>1$, alors $X$, ainsi que la fibre g\'en\'erale de $g$, ne sont pas de type g\'en\'eral (comme on l'a d\'ej\`a remarqu\'e,
les vari\'et\'es de type g\'en\'eral ne poss\`edent pas d'endomorphisme
de degr\'e $>1$). 

4. On a $dim(T)=dim(X)$ si et seulement si $f$ est d'ordre fini.

\medskip

Pour la d\'emonstration, nous avons besoin de quelques pr\'eliminaires.
 
\medskip

{\bf Variantes de la notion d'image par $f$:}

\medskip

Nous allons en utiliser trois. 
Soit $X'\subset X\times X$ le graphe de $f$, et $p,q:X'\to
X$ les projections (de sorte que $f=q\circ p^{-1}$). Ainsi $p$ est une modification propre, et $q$ est
g\'en\'eriquement finie.

1) Image totale: si $A\subset X$, on note $f_*(A):=q\circ p^{-1}(A)$, et
$f^{-1}_*(A):=p\circ q^{-1}(A)$;

2) Image ``usuelle'': $f(A)$ est form\'e des $y\in X$ pour lesquels il
   existe
$x\in X$ tel que $f$ est defini en $x$ et $f(x)=y$;

3) Image ``propre'': pour $A\subset X$ ferm\'e de Zariski , sans composante
 irr\'eductible
dans le lieu d'ind\'etermination de $f$,  $\bar f(A)$ est la
fermeture de Zariski de $f(A)$ dans $X$.

\medskip

On observera que $A\cap f^{-1}_*(B)\neq \emptyset$ si et seulement si
$f_*(A)\cap B\neq \emptyset$.

Si $A$ est un ferm\'e de Zariski (dans $X$), $f_*(A),f^{-1}_*(A)$ le sont aussi. 

Si, de plus, $A\neq X$, alors $f_*(A),f^{-1}_*(A)$ sont aussi des
ferm\'es de Zariski stricts de 
$X$ (car $p,q$ sont g\'en\'eriquement finies).

\medskip

{\bf Lieux d'ind\'etermination:}

\medskip

Soit $I',J'\subset X'$ d\'efinis respectivement par:

$I':=\lbrace x'\in X'\vert p$ non submersive en $x'\rbrace$, et:

$J':=\lbrace x'\in X'\vert p$ ou $q$ non submersive en $x'\rbrace$. 

Alors: $I'\subset J'$.

\

On pose $I:=p(I'), J:=p(J')$: ils sont ferm\'es de Zariski dans $X$,
par le th\'eor\`eme de 
Remmert; $I$ est
le lieu d'ind\'etermination de $f$; et $J$ est le lieu en lequel $f$ n'est
pas simultan\'ement holomorphe et submersive.

\medskip

Posons  $I_{\infty}=\cup _{m\in \Z}(f_*)^mI$, et $J_{\infty}=\cup
_{m\in \Z}(f_*)^mJ$.

Soient $X^+_{\infty}=X-I_{\infty}$, et $ X_{\infty}=X-J_{\infty}$.

Ainsi, $f^m$ est holomorphe et submersive en $x$ pour
tout $m>0$ si $x\in X_{\infty}$;
de plus, $X_{\infty}$ est $f$-invariante.

\medskip

{\bf Adh\'erence de Zariski des f-orbites:}

\medskip

Pour $x\in X_{\infty}$, on
notera $Z^+(x)$ l'adh\'erence de Zariski dans $X$ de la $f$-orbite de
$x$ (constitu\'ee des $f^m(x)$, pour $m\geq 0$). Soit $d(x)=dim(Z^+(x))$.

Aucune composante irr\'eductible de $Z^+(x)$ n'est contenue dans $I_{\infty}$, ni donc dans $I$.

On a alors $$\bar f(Z^+(x))=Z^+(f(x))\subset Z^+(x).$$

On note $Z(x)$ la r\'eunion des composantes irr\'eductibles de
$Z^+(x)$ qui sont de dimension $d(x)$. 

Puisque $f$ est g\'en\'eriquement submersive le
long de la $f$-orbite de $x$ si $x\in X_{\infty}$, on a aussi:
$$\bar f(Z(x))=Z(f(x))\subset Z(x).$$ 

La suite des
$(\bar f)^j(Z^+(x))=Z^+(f^j(x))$ est donc d\'ecroissante pour l'inclusion, et par cons\'equent
stationnaire. Soit $m(x)$ le plus petit des entiers $j\geq 0$ tels que
$$Z^+(f^j(x))=Z^+(f^{j+1}(x))=Z^+(f^{m(x)}(x)).$$ 
Cette \'egalit\'e est donc satisfaite pour tout $j\geq m(x)$.

Autrement dit, $Z^+(y)=Z(y)=\bar f(Z(y))$ si $y=f^j(x)$, pour 
tout $j\geq m(x)$, et tout $x\in X_{\infty}$.

Pour $d,m$ entiers positifs, on notera $X_{\infty}^{d,m}$ l'ensemble des
$x\in X_{\infty}$ tels que $d(x)=d$, et $m(x)=m$.

\

Pour $y\in X_{\infty}^{d,m}$, nous avons donc l'\'egalit\'e que nous 
appellerons ``$(Z=Z^+)$'': 

 $$Z(f^m(y))=Z^+(f^m(y)).$$ 

\medskip

{\bf Cat\'egories de Baire analytiques:}

\medskip

On dira que $E\subset X$ est de {\it premi\`ere cat\'egorie de Baire
analytique dans $X$ (ce qui sera not\'e: $E\in \cal B$$(X)$}) si $E$ n'est pas contenu dans une
r\'eunion d\'enombrable de ferm\'es analytiques d'int\'erieurs vides de
$X$. Nous dirons aussi que $F\subset X$ est {\it negligeable}, si $F\notin \cal B$$(X)$.

\

{\bf Exemples:}

1. Si  $E\notin \cal B$$(X)$, alors $\bar E\in \cal B$$(X)$, o\`u $\bar E$ est le compl\'ementaire de $E$ dans $X$ .

2. Si $E$ est r\'eunion d\'enombrable des $E_m$, et si $E\in \cal B$$(X)$,
alors l'un au moins des $E_m\in \cal B$$(X)$.

3. Donc $X_{\infty}\in \cal B$$(X)$, et l'un au moins des $X_{\infty}^{d,m}\in \cal B$$(X)$.

\

{\bf Lemme 4.2:} {\it Si $A\subset X_{\infty}\in \cal B$$(X)$, 
alors $f^m(A)\in \cal B$$(X)$, pour tout $m\geq 0$.}

\medskip

{\it D\'emonstration:} Sinon, on aurait: 
$f^m(A) \subset \cup_{j=0}^{j=\infty} B_j$, les $B_j$ \'etant des 
ferm\'es de Zariski
stricts dans $X$. Donc $A\subset \cup_{j=0}^{j=\infty} ((f_*)^{-m}(B_j))$. 
Contradiction, puisque $A\in \cal B$$(X)$.

\

{\bf Corollaire 4.3:} {\it Il existe $d,m$ tels que si  
$R=R_{d,m}$ est la r\'eunion
 des $Z(x)$, pour $x\in X_{\infty}^{d,m}$, alors $R\in \cal B$$(X)$.}
 
 \medskip

{\it D\'emonstration:}  En effet, $R$ contient $f^m(X_{\infty}^{d,m})$, et
$f^m(X_{\infty}^{d,m})\in \cal B$$(X)$, par le lemme et l'exemple 3 ci-dessus.

\medskip

{\bf L'espace des cycles:}

\medskip

Soit $C(X)$ l'espace des cycles analytiques compacts de dimension
pure d'un espace analytique 
complexe (\cite{Bar});
lorsque $X$ est projectif, il s'agit simplement de la vari\'et\'e de Chow
de $X$.

Si $X$ est compacte k\" ahl\'erienne, les composantes connexes de $X$ sont
compactes (\cite{L}).

Si $t\in C(X)$, on notera $Z_t$ le cycle analytique compact et de
dimension pure param\'etr\'e par $t$.

Soit $T$ un sous-ensemble analytique compact et irr\'eductible de
$C(X)$ dont le point g\'en\'erique $t$ param\`etre un cycle $Z_t$ dont les
composantes irr\'eductibles sont toutes de multiplicit\'e $1$ et non
contenues dans $I$.

On notera $G=G_T\subset T\times X$ son graphe d'incidence,  $r=r_T:G\to
X$ et $s=s_T:G\to T$ les projections. On notera aussi $L_T:=r_T(G_T)\subset X$
le {\it lieu de $T$}.

On dira que $Z_t$ est $f$-stable si aucune composante de $Z_t$ n'est
contenue dans $I$, et si $\bar f(Z_t)=Z_t$, o\`u $\bar f$ d\'esigne l'image propre de $Z_t$ par $f$.

\smallskip

Par exemple, si $x\in f^m(X^{d,m}_{\infty})$, alors $Z(x)$ est
$f$-stable de dimension $d$.

\

{\bf Lemme 4.4:} {\it Soit $C_{f-stab}(X)$ l'ensemble des $t\in C(X)$ tels que $Z_t$
soit $f$-stable, et  soit $\bar C_{f-stab}(X)$ son adh\'erence de Zariski dans $C(X)$. 

Alors: pour toute composante irr\'eductible $T$ de  $\bar
C_{f-stab}(X)$, 
l'intersection $T_{f-stab}:=T\cap  C_{f-stab}(X)$ contient un ouvert de Zariski (non vide) de $T$.}

\smallskip

(De sorte que $T_{f-stab}$ est  Zariski dense dans $T$ et
le point g\'en\'erique $t$ de $T$ param\`etre un cycle $f$-stable)

\medskip 
 
{\it D\'emonstration:} 
D'abord, il est clair que ``ne pas avoir de composante dans $I$'' et ``\^etre
r\'eduit'' sont des conditions ouvertes (de Zariski).

Soit $G''\subset T\times X$ la fermeture de Zariski de
 $G'=(id_T\times f)(G_T)$.   
 
 Apr\`es aplatissement g\'eom\'etrique et modification de $T$,
on peut supposer que $G''$ est \'equidimensionnel sur $T$, normal. 
Donc $G''$ est le graphe
d'incidence d'une famille analytique de cycles (g\'en\'eriquement images 
stricts de cycles de $G_T$), not\'es $Z_t''$, de $X$, et
param\'etr\' ee par $T$. Soit $U_T$ l'ouvert de Zariski dense de $T$
constitu\'e des $t$ tels que 
$Z_t''$ et $Z_t$ sont r\'eduits et n'ont pas de composante
irr\'eductible dans $I$. 

Alors $T_{f-stab}\cap U_T$ est \'egal \`a $( F^+\cap F^-)$, $F^+$
(resp. $F^-$) 
\'etant constitu\'e des
 $t$ pour lesquels $Z_t\subset Z_t''$ (resp. $Z_t''\subset
 Z_t$). Puisque $F^+$ et $F^-$ sont des ferm\' es de Zariski dans $T$, le lemme est d\'emontr\'e.
 
 \

 {\bf Corollaire 4.5:} {\it Il existe une famille d\'enombrable de sous-ensembles
analytiques compacts et irr\'eductibles $T_k,k\geq 0$ de $C(X)$ tels que,
pour chaque $k$, le membre g\'en\'erique $Z_t$ de $T_k$ soit $f$-stable,
et tels que tout cycle $f$-stable de dimension pure de $X$ soit
param\'etr\'e par un point de l'un des $T_k$.}

 \

Pour chaque $k\geq 0$, soit $L_k$ le lieu de la famille $T_k$. Rappelons
que c'est un ferm\'e de Zariski de $X$.

On dira que $T_k$ est {\it couvrante} si $L_k=X$. Il existe des $T_k$
couvrantes (par exemple: celle r\'eduite \`a $Z_t:=X$.

\smallskip

On utilise la notation $R_{d,m}$ introduite dans le corollaire 4.3.

\

 {\bf Lemme 4.6:} {\it Pour chaque couple d'entiers $(d,m)$ tel que
$R_{d,m}\in \cal B$$(X)$,
il existe un entier $k\geq 0$ v\'erifiant les propri\'et\'es
suivantes:

1. $T_k$ est couvrante;

2. Il existe $E_k^{d,m}\in \cal B$$(X)$, 
$E_k^{d,m}\subset f^m(X^{d,m}_{\infty})$,
tel que $Z(x)=Z_t$, 
pour un certain $t=t(x)\in T_k$, d\`es que $x\in E_k^{d,m}$.}

\medskip

{\it D\'emonstration:} Soit $L\subset X$ la r\'eunion d\'enombrable des
$L_k$ pour lesquels $T_k$ n'est pas couvrante. On choisit $(d,m)$ tel que
$R_{d,m}\in \cal B$$(X)$. Il existe un tel couple, par le corollaire 4.3.
 Soit $\tilde{E}^{d,m}=(f^m(X^{d,m}_{\infty}))\cap
(R_{d,m}-L)$. C'est un ensemble de premi\`ere cat\'egorie de Baire analytique 
dans $X$, dont le compl\'ementaire est 
negligeable dans $f^m(X^{d,m}_{\infty})$. 
Alors, si $x\in \tilde{E}^{d,m}$,  $Z(x)$ est un cycle
de dimension pure $d$ qui est $f$-stable, donc param\'etr\'e par un 
$t\in T_{k(x)}$, pour un
certain $k(x)$, et on a donc: $Z(x)=Z_t$. Si l'on note $\tilde{E}^{d,m}_k$ 
l'ensemble des $x\in \tilde{E}^{d,m}$ pour lesquels $t(x)\in T_k$, il existe 
(par d\'enombrabilit\'e de l'ensemble des $k$) au moins un $k$ pour 
lequel $\tilde{E}^{d,m}_k\in \cal B$$(X)$. 
On pose $E_k^{d,m}=\tilde{E}^{d,m}_k$ pour un tel $k$. 
La famille $T_k$ est couvrante,
puisque $x\in X-L$, par hypoth\`ese, et $E_k^{d,m}$ 
satisfait les propri\'et\'es annonc\'ees.

\medskip

Nous pouvons maintenant d\'emontrer le th\'eor\`eme annonc\'e au
d\'ebut
de cette partie.

\

{\it Preuve du th\'eor\`eme 4.1:}

\smallskip

On choisit d\'esormais $k,d,m, E$ comme ci-dessus. On pose $T:=T_k$, et on
note $r:G\to X$, $s:G\to T$ les projections , si $G\subset T\times X$ est
le graphe d'incidence de la famille $T=T_k$.

D\'emontrons les assertions suivantes:

\

1. $r$ est une modification propre. 

Soit $g=s\circ r^{-1}:X\dasharrow T$: c'est
une application m\'eromorphe surjective de fibre g\'en\'erique $Z_t$.

2. $g\circ f=g$.

3. Pour $x\in X$ g\'en\'eral,  $Z(x)=Z_{g(x)}$.

\

Montrons l'assertion 1. La surjectivit\'e de $r$
r\'esulte de ce que $T$ est couvrante. De plus, pour $t\in T$
g\'en\'erique, 
$Z_t$ n'a aucune composante irr\'eductible contenue dans $I$. 
On en d\'eduit que si $r$ n'est pas bim\'eromorphe, et si $x\in X$
est g\'en\'erique (c.\`a.d. en dehors d'un ferm\'e de Zariski
strict), 
il existe (au moins) deux $t\neq t' \in T$ tels que $x\in Z_t\cap
Z_{t'}$, avec 
$Z_t$ et $Z_{t'}$ sans composante irr\'eductible contenue dans $I$.

Montrons maintenant que si $x\in E$, $E\in \cal B$$(X)$ ad\'equat,
alors la fibre de 
$r$ au-dessus de $x$ ne peut avoir deux tels points $t,t'$. 
Cette propri\'et\'e entra\^{\i}nera bien que $r$ est une modification propre.

Nous choisissons $x\in E^{d,m}_k \subset f^m(X_{\infty}^{d,m})$, pour
$(d,m,k)$ et $E^{d,m}_k$ comme dans 
le lemme pr\'ec\'edent, et raisonnons par l'absurde.

Soient $t\neq t'$ deux points distincts de $T$, tels donc que: $x\in
Z_t\cap Z_{t'}\neq Z_t$. 
On peut supposer,  par notre choix de $x$, que $Z_t=Z(x)$.
Puisque $x=f^m(y),y\in
X_{\infty}^{d,m}$,
on a:

$Z_t=Z(x)=Z(f^m(y))=Z^+(f^m(y))=Z^+(x)$, ceci par l'\'egalit\'e $(Z=Z^+)$. 

Donc $Z_t \cap Z_{t'}\neq Z_t$ contient $x$, et est $f$-stable; autrement
dit, $Z_t \cap Z_{t'}\neq Z_t$ contient $Z^+(x)$.

Mais ceci contredit l'\'egalit\'e $Z_t=Z^+(x)$.

Donc $r$ est bien une modification propre.

\

La propri\'et\'e 2 r\'esulte de ce que $Z_t$ est $f$-stable, pour $t\in T$
g\'en\'erique.

\

Montrons la propri\'et\'e 3:  si elle n'est pas satisfaite, il existe
un sous-ensemble $U\in {\cal B}(X)$, tel que pour $x\in U$, 
$Z(x)\subsetneq G(x)$, o\`u $G(x)$ d\'enote (l'unique) fibre de $g$ passant
par $x$. Soit $U_{\infty}=U\cap X_{\infty}$, et $U^{d,m}_{\infty}=
U\cap X^{d,m}_{\infty}$. Il existe un couple $m',d'$ tel que 
$U^{d',m'}_{\infty}\in {\cal B}(X)$, et donc 
$f^m(U^{d',m'}_{\infty})\in {\cal B}(X)$. De plus, pour $y$ dans un sous-
ensemble ad\'equat (et dans ${\cal B}(X)$) $E'$ de $f^m(U^{d',m'}_{\infty})$,
nous avons toujours $Z(y)\subsetneq G(y)$.

En rep\'etant l'argument du lemme 4.6 et de la preuve de la propri\'et\'e 1,
nous obtenons une seconde fibration $f$-invariante $g':X\dasharrow T'$,
telle que la fibre $G'(y)$ par $y\in E'$ coincide avec $Z(y)$.
Au m\^eme temps, la fibre $G'(x)$ par $x\in E=E^{d,m}_k$ contient
$Z^+(x)\supset Z(x)$, par $f$-invariance de $g'$. Ainsi, 
pour les sous-ensembles $E,E'\in \cal B$$(X)$, on
a: $Z(x)=Z_{g(x)}\subset Z_{g'(x)}=G'(x)$, pour $x\in E$, et
$Z(y)=Z_{g'(y)}\subset Z_{g(y)}=G(y)$, pour $y\in E'$.  Puisque ces
inclusions d\'efinissent des sous-ensembles analytiques 
ferm\'es de $T,T'$, on a
ainsi une contradiction si $T\neq T'$ (en tant qu'une famille de cycles
sur $X$).

\

Nous avons donc montr\'e que la fibre $G(x)$ de $g$ passant par 
$x\in X$ g\'en\'eral
 est \'egale \`a $Z(x)$. Pour finir, remarquons que cette fibre
contient $Z^+(x)$ par $f$-invariance; d'o\`u $G(x)=Z(x)=Z^+(x)$.
Le th\'eor\`eme est d\'emontr\'e.

\

Pour conclure cette partie, consid\'erons \`a 
nouveau l'exemple de Voisin: soit $X={\cal F}(V)$ la
vari\'et\'e des droites de $V$, cubique lisse dans $\p^5$.
D'apr\`es \cite{BD}, nous avons un isomorphisme d'Abel-Jacobi 
$$H^4(V,\Z)\cong H^2(X,\Z),$$
induisant un isomorphisme de structures de Hodge. D'autre part,
un r\'esultat de Terasoma (\cite{T}) garantit l'existence
d'une cubique lisse $V$ d\'efinie sur un corps de nombres, et telle que
$H^{2,2}(V)\cap H^4(V,\Z) =\Z$. Ceci entra\^{\i}ne \' evidemment la

\

{\bf Proposition 4.7} {\it Il existe un $X$ comme ci-dessus,
  d\'efini
sur un corps de nombres et tel que $Pic(X)\cong \Z$.}

\

La $f$-orbite d'un point g\'en\'eral (sur $\C$)  d'un tel $X$ est 
Zariski-dense par les r\'esultats de cet article. Si l'on pouvait
obtenir un \'enonc\'e semblable sur $\bar \Q$, on aurait la densit\'e
potentielle
de points rationnels sur $X$. Mais il est clair que notre argument
ne fonctionne pas sur $\bar \Q$, puisque nous consid\'erons comme
``n\'egligeable'' toute r\'eunion d\'enombrable de ferm\'es stricts,
et $\bar \Q$ lui-m\^eme est d\'enombrable. Il serait int\'eressant
de trouver des moyens de distinguer les r\'eunions d\'enombrables de
ferm\'es
``arbitraires'' des r\'eunions d\'enombrables ``naturelles'' de
ferm\'es, obtenues par exemple par une construction g\'eom\'etrique
d\'efinie sur un corps de nombres \`a partir d'un seul ferm\'e.

\

Pour un corps non-d\'enombrable $F$, $car(F)=0$, la situation est diff\'erente.
En effet, il est clair que, dans le Th\'eor\`eme 4.1, on peut
prendre pour $X$ une vari\'et\'e projective (irr\'eductible) 
sur $\bar F$. On en d\'eduit facilement de nouveaux exemples de la densit\'e
potentielle sur un corps de fonctions (cette observation est
motiv\'ee par le texte \cite{HT2}):

\

{\bf Proposition 4.8} {\it Soit $F$ un corps non-d\'enombrable (par
 exemple, $F=\C(C)$, o\`u $C$ est une courbe complexe).
 Soit $X$ une vari\'et\'e lisse projective
sur $\bar F$, v\'erifiant $K_X=0$, $Pic(X)=\Z$, et munie d'un endomorphisme rationnel
$f:X\dasharrow X$ de degr\'e $d>1$. Alors il existe une extension
  finie
$L$ de $F$, telle que l'ensemble $X(L)$ est Zariski-dense dans $X$.}

\

En effet, les conditions $K_X=0$, $Pic(X)=\Z$ assurent par 2.1 que $f$ ne
peut
pas pr\'eserver une fibration non-triviale. 
Par 4.1, il existe donc une r\'eunion d\'enombrable $R$ de ferm\'es de Zariski stricts de $X$ tel que tout point
$x\in X$, $x\notin R$ est de $f$-orbite d\'efinie et Zariski-dense. 
Puisque $F$ n'est
pas d\'enombrable, $R$ ne contient pas $X(\bar{F})$. 
Il existe donc un point
$x\in X$ dont l'orbite iter\'ee est Zariski-dense, et il suffit
de prendre pour $L$ une extension finie de $F$ 
sur laquelle $X$, $x$ et $f$ sont
d\'efinis.

\

Les exemples de Voisin fournissent tout de suite des vari\'et\'es
$X$ sur $\bar F$, $F=\C(C)$, non-isotriviales en tant que familles
sur $C$, 
satisfaisant
les conditions de la proposition 4.8. Il suffit de prendre 
$C\subset \p(H^0(\p^5,{\cal O}(3)))$ assez g\'en\'erale, et consid\'erer
la famille ${\cal X}\rightarrow C$, dont chaque fibre $X_c$ est la 
vari\'et\'e de Fano des droites de la cubique $V_c$ correspondant
\`a
$c\in C$. On prend
ensuite
pour $X$ la fibre g\'en\'erique g\'eom\'etrique $X_{\bar F}$. 

Alternativement, on peut prendre $x_c\in X_c$ ($c\in C$ assez
g\'en\'eral), tel que la $f_c$-orbite de $x_c$ est Zariski-dense
dans $X_c$ (ici, $f_c$ d\'enote l'endomorphisme de Voisin de $X_c$).
On choisit une multisection $\tilde{C}$ passant par $x_c$; apr\`es
un changement de base $\tilde{C}\rightarrow C$, les sections de
${\cal X}$ sont Zariski-denses, c'est-\`a-dire, les
$\C(\tilde{C})$-points de $X_{\bar F}$ sont Zariski-denses.

\

\section{Base d'une fibration hyperk\" ahl\'erienne}

Soit $g:X\dasharrow B$ une fibration m\'eromorphe, avec $X$  
hyperk\"ahl\'erienne irr\'eductible (non n\'ecessairement
projective) de dimension complexe $2n$, et $B$ lisse k\"ahlerienne. 
On suppose $g$ non triviale: $0<dim(B)<2n$. Comme auparavant,
on d\'esigne par $\sigma$ la forme symplectique sur $X$.

\medskip

{\bf Question:} {\it $B$ est-elle alors une vari\'et\'e rationnelle?}

\medskip

{\bf{Remarques:}}

\medskip

1. $B$ est une vari\'et\'e projective. En effet, sinon elle admet une
$2$-forme holomorphe non-nulle (consequence du crit\`ere de projectivit\'e 
de Kodaira et de la d\'ecomposition de Hodge). Son image
r\'eciproque $u$ par $g$ est non-nulle et n'est pas de rang maximal. Ceci
contredit le fait que $X$ est irr\'eductible (et donc que $u$ et $\sigma$ sont
proportionnelles).

2. Supposons que $g$ ne soit pas presque holomorphe. Alors $B$ est
unir\'egl\'ee. Si donc $B$ n'est pas rationnellement connexe, et si 
$r:B\dasharrow R$ est son quotient rationnel (voir \cite{C2}, \cite{KMM}), 
alors $R$ n'est pas unir\'egl\'ee par \cite{GHS}, et
$h=rg:X\dasharrow R$ est une fibration non-triviale presque holomorphe.
Ainsi, ou bien $B$ est rationnellement connexe, ou bien il existe sur $X$
une fibration $h:X\dasharrow R$ non triviale et presque holomorphe.

3. Si $B$ est rationnelle, $B$ est bien s\^ur rationnellement connexe, et
les tenseurs holomorphes covariants sur $B$ s'annulent tous.
Conjecturalement, cette
propri\'et\' e entra\^{\i}ne  la connexit\'e rationnelle de $B$.

\medskip

Cette propri\'et\'e d'annulation est effectivement verifi\'ee, par
une g\'en\'eralisation de 
l'argument tr\`es simple de la remarque 1:

\

{\bf Proposition 5.1} 
{\it Soit $g:X\dasharrow B$ comme ci-dessus. Soient $m>0,p>0$ des
entiers.

Alors $H^0(B,Sym^m(\Omega^p_B))=0$.}

\

La proposition r\'esulte du lemme suivant:

\

{\bf Lemme 5.2} {\it Soit $X$ hyperkahl\'erienne irr\'eductible.
Pour $m,p$ des entiers strictement positifs, l'application $c:Sym^m(H^0(X,\Omega^p_X))\to  H^0(X,Sym^m(\Omega^p_X))$ est un isomorphisme. 

Le groupe $H^0(X,Sym^m(\Omega^p_X))$ est donc nul si 
$p$ est impair, et est \'egal \`a
$\C(\sigma^{\wedge
q})^{\otimes m}$ si $p=2q$ est pair.}

\

Ce lemme para\^{\i}t bien connu, r\'esultant de l'existence d'une  m\'etrique 
de K\" ahler Ricci-plate sur $X$ et du principe de Bochner
(voir par exemple \cite{Kob} pour un raisonnement semblable).

\

Supposons maintenant qu'il existe
 $s\in H^0(B,Sym^m(\Omega^p_B))$ non-nul. Soit
$0\neq s':=g^*(s)\in H^0(X,Sym^m(\Omega^p_X))$: ce tenseur est bien d\'efini
par le th\'eor\`eme
de Hartogs.
Le lemme 5.2 montre que $p=2q$ est pair, et qu'il existe
$\lambda \in \C$, non nul, tel que 
$s'=\lambda (\sigma^{\wedge q})^{\otimes m}$. 
Mais un calcul standard en coordonn\'ees locales montre que ceci est
impossible d\`es que $dim(B)<2n$: au point g\'en\'erique $x\in X$, 
on construit facilement un \'el\'ement $\tau\in (\Lambda^{2q}T_x
 X)^{\otimes m}$ annul\'e par tous les $g^*(s)$ mais non par $(\sigma^{\wedge q})^{\otimes m}$.

\

Si $dim(B)=2$, cette annulation implique que $B$ est rationnelle,
par le crit\`ere de Castelnuovo. En g\'en\'eral, on dit que 
$\kappa_+(B)=-\infty$ si
$\kappa(B')=-\infty$ pour tout $B'$ domin\'e par $B$ (c.\`a.d. tel qu'existe
$h:B\dasharrow B'$ m\'eromorphe dominante). Le th\'eor\`eme 5.1 implique, bien s\^ur,
que dans notre situation, $\kappa_+(B)=-\infty$. 

De l'existence du quotient rationnel 
(\cite{C2},\cite{KMM}) et de \cite{GHS},
on d\'eduit que $\kappa_+(B)=-\infty$ si et seulement si $B$ est
rationnellement connexe, pourvu que les vari\'et\'es de dimension
$\leq dim(B)$ avec $\kappa=-\infty$
soient unir\'egl\'ees. Comme cette assertion est demontr\'ee en dimension $3$, 
nous avons:

\

{\bf{Corollaire 5.3}} {\it Si $g:X\dasharrow B$ est comme ci-dessus, 
$B$ est rationnelle si
$dim(B)=2$, et $B$ est rationnellement connexe si $dim(B)=3$.}

\section{Endomorphismes et classification}

Soit $X$ une vari\'et\'e K\" ahl\'erienne compacte lisse et connexe,
et soit 
$f:X\dasharrow X$ un endomorphisme m\'eromorphe dominant de $X$ de
degr\'e $d\geq 1$. 
On se propose ici de montrer que $f$ respecte certaines fibrations 
intrins\`equement attach\'ees \`a $X$.

Parmi celles-ci, figurent le quotient rationnel $r:X\dasharrow R$
(\cite{C2},\cite{KMM}), la 
fibration d'Iitaka-Moishezon $\phi:X\dasharrow J$ (si
$\kappa(X)\geq 0$) , 
et le {\it coeur} $c:X\dasharrow C$ (\cite{C3},3.1,3.3, p. 544; 5.8,
p. 579). Ils sont d\'efinis \`a \'equivalence bim\'eromorphe pr\`es
seulement. Nous modifierons donc les mod\'eles choisis au gr\'e des besoins.

\

{\bf Proposition 6.0} {\it Dans la situation pr\'ec\'edente, 
il existe des endomorphismes dominants $r_f:R\to R$ et $\phi_f:J\to
J$, tels que 
$r\circ f=r_f\circ r$ et $\phi\circ f=\phi_f\circ r$ .

De plus, $\phi_f:J\to J$ est bim\'eromorphe.}

\

{\bf D\'emonstration:} Soit $F$ une fibre g\'en\'erale de $r$. Alors
$F$ est rationnellement connexe. Donc $f(F)$ l'est aussi; de plus,
$dim(F)=dim(f(F))$. De la propri\'et\'e universelle  du quotient
rationnel r\'esulte donc que $f(F)$ est une fibre de $r$.

Consid\'erons maintenant la fibration $\phi$. L'application
$f$ induit, par image r\'eciproque, 
un automorphisme de $H^0(X,mK_X)$ pour tout $m$ entier. En effet,
l'image  r\'eciproque envoit injectivement, et donc surjectivement, $H^0(X,mK_X)$ dans lui-m\^eme.
Nous avons alors, pour un certain $m$, un diagramme
commutatif

$$
\begin{array}{ccc}
{X} & \stackrel{\phi}{\rightarrow} & {\p(H^0(X,mK_X)^*)} \\
f{\downarrow} &&  \downarrow \\
{X} & \stackrel{\phi}{\rightarrow} & {\p(H^0(X,mK_X)^*)}
\end{array}
$$
 
L'application $\phi_f$ est la restriction sur $J=Im(\phi)$ de la seconde fl\`eche verticale
(l'automorphisme correspondant de $\p(H^0(X,mK_X)^*)$). Elle est
\'evidemment bim\'eromorphe.

\

{\bf Remarque:} La structure de $r_f$ peut \^etre arbitrairement
compliqu\'ee comme des exemples tr\`es simples (produits) 
le montrent. 

\

Pour le ``coeur", qui ``scinde" naturellement 
la structure de $X$ en
ses deux composantes antith\'etiques (``sp\'eciale" et `` de type
g\'en\'eral"), 
la situation est tr\`es claire:

\

{\bf Th\'eor\`eme 6.1}   {\it  Soit $X$ une vari\'et\'e 
k\"ahl\'erienne compacte 
lisse et connexe, et $c:X\to C$ le ``coeur" de $X$ . Si $f:X\dasharrow X$ est un endomorphisme m\'eromorphe dominant de $X$ de degr\'e $d\geq 1$, alors:

1. Il existe $c_f:C\dasharrow C$ bim\'eromorphe telle que $c\circ f=c_f\circ c$;

2. Il existe $N>0$ entier tel que $f^N$ pr\'eserve $c$ (ie: $c\circ f^N=c$).

En particulier:

3. $c$ se factorise par l'application $g_N:X\dasharrow T_N$ du th\'eor\`eme 4.1
correspondant \`a une certaine puissance $f^N$ de $f$;

4. Si l'application $g_N:X\dasharrow T_N$ est constante (autrement dit, si un
point $x$ de $X$ 
a une f-orbite Zariski dense dans $X$), alors $X$ est ``sp\'eciale"
au 
sens de \cite{C3}, 2.1, p.527.}

\

{\bf D\'emonstration:} Les propri\'et\'es 3 et 4 sont des
cons\'equences imm\'ediates 
des propri\'et\'es 1 et 2 que nous montrons maintenant.

Rappelons que les fibres g\'en\'erales $S$ de $c$ sont sp\'eciales
et sa base orbifolde de 
type g\'en\'eral (propri\'et\'es qui d\'eterminent $c$). Comme l'image de $S$,
fibre g\'en\'erale de  
$c$, par $f$ est sp\'eciale et de m\^eme dimension que $F$, $f(S)$
doit \^etre une fibre de $c$ (par construction de $c$). 
D'o\`u l'existence de $c_f$ telle que $c\circ f=c_f\circ c$. 

Il reste \`a voir que $c_f$ est bim\'eromorphe et d'ordre fini 
(c.\`a.d. que $(c_f)^N=id_C$, pour un $N>0$ ad\'equat).

Pour la premi\`ere assertion ($c_f$ bim\'eromorphe), rappelons que
si $(C/\Delta(c))$ est la base orbifolde de $c$, et si
$p:=dim(C)>0$, alors $L_X:=\lceil c^*(K_C+\Delta(c))\rceil\subset
\Omega_X^p$ est (sur un mod\`ele bim\'eromorphe ad\'equat)
localement libre de rang $1$, et de dimension de Moishezon-Iitaka
\'egale \`a $p$ (autrement dit: un faisceau {\it de
  Bogomolov}). Voir \cite{C3}, Section 2.6 pour d\'etails. 

Un tel faisceau $L_X$ est unique, puisque l'on a unicit\'e du ``coeur"
$c$, 
et que $L_X$ est d\'etermin\'e par $c$. Donc $f^*(mL_X)\subset
mL_X$, pour tout entier $m>0$. Autrement dit, $f^*$ induit un
automorphisme de $H^0(X, mL_X)$. Puisque, pour $m$ assez  grand et
divisible, le syst\`eme lin\'eaire $|mL_X|$ (disons, de dimension
$M$) 
fournit $c$, nous pouvons
conclure,
par le m\^eme diagramme que dans la preuve du 6.0, qu'il existe 
$g\in PGl(M+1,\C)$ avec $g\circ c=c\circ f$. L'application $c_f$,
bim\'eromorphe,
est la restriction de $g$ sur l'image $C$ de $c$.

\

Pour d\'emontrer la finitude, \'etablissons les lemmes suivants
(dont le second est une g\'en\'eralisation orbifolde directe 
de (\cite{U}), 14.3), montrant que le groupe des automorphismes
bim\'eromorphes 
d'une vari\'et\'e de type g\'en\'eral est fini:

\

{\bf Lemme 6.2}  {\it Dans la situation pr\'ec\'edente, $g_*(\Delta(c))=\Delta(c)$. (Autrement dit: l'automorphisme $g$ pr\'eserve la base orbifolde de $c$)}

\

{\bf D\'emonstration:} Soit $u:X'\to X$ une modification, avec $X'$
lisse, 
telle que $f':=f\circ u:X'\to X$ soit holomorphe. 

Nous avons: $g_*(\Delta(c))=\Delta(g\circ c)=\Delta(g\circ c\circ
u)=\Delta(c\circ f')\geq \Delta(c)$, 
o\`u, pour des $\Q$-diviseurs $A,B$, $A\geq B$ signifie que $A-B$
est effectif ou nul. 
Puisque $g$ est un automorphisme de $C$, ceci entra\^{\i}ne que
$g_*(\Delta(c))=\Delta(c)$.

\

{\bf Lemme 6.3}  {\it Soit $c:X\to C\subset \p^M$ le mod\`ele
  pr\'ec\'edent du coeur  
de $X$, fourni par le syst\`eme lin\'eaire $mL_X$, pour $m>0$ assez
  grand et divisible. 
Le groupe $G_X$ des $g\in \p GL(M+1,\C)$ pr\'eservant $C$, 
tels que $g_*(\Delta( c))=\Delta(c)$, est fini.}

\

{\bf D\'emonstration:} Soit $H_C$ le sous-groupe des $g\in \p GL(M+1,\C)$ 
pr\'eservant $C$. Le groupe  $G_X$ s'injecte naturellement dans $H_C$. 
Le
th\'eor\`eme 14.1 de \cite {U} 
(bas\'e sur un r\'esultat de Rosenlicht) montre que si $H_C$ n'est
pas fini, 
il contient un sous-groupe lin\'eaire alg\'ebrique $K\subset \p
GL(2,\C)$ de dimension $1$, 
et que $C$ est bim\'eromorphe \`a $W\times \p^1$, pour une certaine
vari\'et\'e $W$. 
Ceci de telle sorte que les sous-vari\'et\'es $\lbrace w\rbrace\times \p^1$ soient les 
adh\'erences de Zariski des orbites de $K$ agissant sur
$C$. Changeant de mod\`ele bim\'eromorphe
 pour $X$ et $c$, nous supposerons que $c':X\to C'=W\times \p^1$ est
 le coeur de $X$.
 Le diviseur canonique orbifolde $K_{C'}+\Delta(c')$ de la base
 orbifolde de $c'$ est 
donc encore de type g\'en\'eral (puisque le coeur est une fibration
de type g\'en\'eral, 
ce qui signifie que les bases orbifoldes de {\it tous} ses mod\`eles
bim\'eromorphes sont de 
type g\'en\'eral), de sorte que pour $w\in W$ g\'en\'erique,
$\lbrace w\rbrace\times \p^1$ 
rencontre le support de $\Delta(c)$ en au moins $3$ points. Par
ailleurs l'action naturelle de 
$K\subset \p GL(2,\C)$ sur $\lbrace w\rbrace\times \p^1$ doit
pr\'eserver ce support, 
d'apr\`es 6.2. On a donc une contradiction, et ainsi la finitude  de $G_X$.

\

Les deux lemmes entra\^{\i}nent que $c_f$ est d'ordre fini, C.Q.F.D.

\

Finalement, consid\'erons la fibration d'Iitaka $\phi$.
Nous avons d\'ej\`a remarqu\'e que $\phi_f$ est toujours
bim\'eromorphe. Il est alors naturel de poser la question
suivante: est-il vrai que $\phi_f$ est toujours d'ordre fini?
Bien que, a priori, il semble peu probable que la r\'eponse soit
positive, la question m\'erite l'\'etude. Pour mieux comprendre
la situation, montrons que la r\'eponse est positive pour les surfaces:

\

{\bf Proposition 6.4:} {\it Soit $X$ une surface k\" ahl\'erienne 
 compacte avec $\kappa(X)=1$. Soit $\phi:X\to B$ sa fibration 
d'Iitaka-Moishezon, et $f:X\dasharrow X$ un endomorphisme m\'eromorphe de
 degre 
$d\geq 1$ de $X$. Soit $g:B\to B$ l'automorphisme de $B$ induit par
 $f$ . 
Alors $g$ est d'ordre fini.}

\

{\it Preuve:}  Tout d'abord, l'assertion est claire
si le genre de $B$ est au moins $2$. Si $B$ est elliptique, $g$ est
encore 
d'ordre fini, puisqu'induit par une action lin\'eaire dans un
plongement. Nous pouvons donc supposer que $B=\p^1$, et que $X$ est
relativement minimale.

Soit $M\subset B$ le lieu de fibres
multiples de $\phi$; il est conserv\'e par $g$. 
Montrons que le
lieu de fibres singuli\`eres non multiples 
$S\subset B$ est \'egalement conserv\'e.
Soit $b\in B$ tel que $X_{g(b)}$ soit lisse, donc elliptique. 
Soit $\mu:X'\to X$ une modification telle que $f':=f\circ \mu:X'\to X$
soit holomorphe. Les courbes exceptionnelles de $\mu$ au-dessus de
$X_b$ sont rationnelles, donc envoy\'ees sur des points par
$f'$. D'o\`u 
$f$ est holomorphe au voisinage de $X_b$, et l'une des composantes 
irr\'eductibles de $X_b$ est elliptique. Ceci entraine bien 
(voir par exemple \cite{BPV}, V.7, p. 150) que $(X_b)_{red}$ est
lisse elliptique.

Le lieu $E=M\cup S$ de fibres singuli\`eres de $\phi$ est donc
conserv\'e par $g$; si $g$ est d'ordre infini, $E$ a au plus
deux \'el\'ements. Il est bien connu qu'une fibration en courbes
de base $\p^1$ qui n'a que deux fibres singuli\`eres, est \`a
modules constants: en effet, le rev\^etement universel de $\p^1$ priv\'e de
deux points est $\C$, qui n'admet pas d'application holomorphe non-constante
dans le domaine de p\'eriodes (l'espace de Siegel, ou le demi-plan
sup\'erieur dans le cas d'une fibration elliptique). 

Si $g$ est d'ordre infini, $X$ est donc \`a modules
constants. Le lemme suivant montre que ceci est impossible pourvu que
$\kappa(X)=1$.

\

{\bf Lemme 6.5} {\it Soit $h:X\to B$ une surface k\"ahl\'erienne elliptique 
relativement
minimale de base $B\cong \p^1$, \`a modules constants, 
ayant au plus 2 fibres singuli\`eres. Alors
 $\kappa(X)=-\infty.$}

\

{\it D\'emonstration}: La monodromie de $h$ pr\'eserve une classe de
K\"ahler 
sur
 les fibres lisses. Cette monodromie est donc finie cyclique. Quitte
 \`a faire un
 changement de base fini $C\rightarrow B$, 
ramifi\'e seulement au-dessus de 2 points au
 plus, on peut donc supposer cette monodromie triviale; on a encore
 $C\cong P^1$, et la surface $X'$, obtenue par ce changement de base,
a au plus deux fibres singuli\`eres.

 D'apr\`es Kodaira (\cite{K}, table I,p.604), la famille \'etant
isotriviale de monodromie triviale et relativement minimale, les fibres
singuli\`eres sont donc multiples elliptiques. On d\'eduit alors de la 
formule pour le fibr\'e canonique d'une telle
fibration (\cite{BPV}, 12.1,12.3) que $\kappa(X')=\kappa(X)= -\infty$.

\

{\bf Remarque:} Dans le cas o\`u $X$ n'est pas \`a modules constants,
on peut borner effectivement l'ordre de $g$ en utilisant les courbes
modulaires, et ceci m\^eme sans supposer $\kappa(X)=1$.  
Plus pr\'ecis\'ement: l'ordre de $g$ est major\'e par $210 \nu$, 
o\` u $\nu$ est le degr\'e g\'eom\'etrique de l'application
modulaire 
$j:B\to \Bbb P_1$ associant \`a $b\in B$ g\'en\'erique  l'invariant
modulaire 
$j(X_b)$ de la courbe elliptique $X_b$.

 Voici l'esquisse 
de la d\'emonstration:

Soit $h:X\rightarrow B$ une surface elliptique \`a modules non constants,
$j:B\rightarrow \p^1$ l'application modulaire correspondante, $B^0$ la
partie de $B$ param\'etrant les fibres lisses, $X^0=h^{-1}B^0$. Pour $b\in B$
g\'en\'erique, soit
$G=G_b(f)$ le groupe des automorphismes de $X_b$ 
pr\'eservant les fibres de $f$. C'est un sous-groupe fini de translations 
de $X_b$. 

Il existe un entier $m\geq 1$, tel que $G$
contient toute translation dont l'ordre divise $m$, autrement dit,
le sous-groupe complet de $m$-torsion $T(m)$ du groupe des
translations de $X_b$. Des arguments \'el\'ementaires sur la
structure des groupes ab\'eliens finis montrent que le
quotient $G_b(f)/T(m)$ est cyclique. On pose $T=T(m)$ pour $m$ maximal,
$C=G_b(f)/T(m)=G/T$ le groupe cyclique quotient. 

En consid\'erant la correspondance $\Sigma_f\subset X\times_B X$, form\'ee
des $(x,y)$ tels que $f(x)=f(y)$, on voit facilement que
l'application
quotient par $T$ se prolonge en une application ``de multiplication
par $m$''
$u:X\dasharrow X'=Pic^m(X)$ au-dessus de $B$. 

Ainsi, $f$ se factorise en produit $c\cdot u$, avec $c:X'\dasharrow X$
une application telle que, 
pour la fibre g\'en\'erique fix\'ee $X'_b$, le groupe
des translations de $X'_b$ commutant avec la restriction de $c$ est cyclique
(d'ordre $|C|$). On notera $H$ ce sous-groupe.

L'application $c$ induit $g$ sur la base $B$ et donne un rel\`evement de $j$
en une application $j(c): B\rightarrow X_0(|C|)$. Donc $X_0(|C|)$
est rationnelle. 
Remarquons que $X_0(N)$
est rationnelle pour $15$ valeurs de $N$ seulement (plus 
pr\'ecisement,
on doit avoir
$N\in F$, avec $F=\{1,2,3,4,5,6,7,8,9,10,12,13,16,18,25\}$).

Enfin, l'application $f$ n'est pas, en g\'en\'eral, d\'etermin\'e
par 
le sous-groupe $G$; mais si le
m\^eme $G$ correspond \`a $f_1, f_2:X\dasharrow X$, alors pour les
automorphismes de la base $g_1, g_2: B\rightarrow B$ induits,
$g_1g_2^{-1}$ commute avec $j$. Il en va de m\^eme pour $H$
et $c$.

Maintenant, it\'erons $f$: on a $f^n=c_nu_n$; soit $H_n$ le sous-groupe 
cyclique
correspon- dant. Puisque $B\cong \p^1$, on n'a que 15 possibilit\'es pour
$|H_n|$ - sa valeur doit se trouver dans la liste $F$; et pour chaque
valeur $N$ dans $F$, on n'a que $\psi(N)=N\Pi_{p|N}(1+1/p)$ choix de $H_n$
possibles, puisque $\psi(N)$ est le nombre de sous-groupes cycliques
d'ordre $N$ dans $Aut(X'_b)$. 

Soit $A=\sum_{N\in F}\psi_N=210$; alors $H_n=H_m$ pour certains
$0\leq m,n\leq A$. Puisque $c_i$ induit $g^i$ sur la base $B$,
il s'ensuit que $g^{m-n}$ commute avec $j$. Donc l'ordre de
$g^{m-n}$ est au plus $\nu=deg(j)$, et l'ordre de $g$ 
est au plus $A\nu=210\nu$.

\

{\bf Remarque:} La ``r\'eciproque'' de 6.1.4 est fausse pour les
surfaces elliptiques de base
$\p^1$,
de dimension canonique 1, admettant une section (donc sp\'eciales). Nous
avons d\'ej\`a observ\'e 
que, pour
une telle surface, $\phi_f$ est d'ordre fini; $f$ ne peut donc pas
avoir
une orbite it\'er\'ee Zariski-dense.

\

{}

\end{document}